\documentclass[]{article}

\usepackage[T1]{fontenc}
\usepackage[latin1]{inputenc}
\usepackage[english]{babel}
\usepackage{latexsym}
\usepackage{latexsym,amssymb,amsmath,epsfig, graphicx, epstopdf}

\newtheorem{theorem}{Theorem}[section]

\newtheorem{definition}{Definition}[section]
\newtheorem{proposition}{Proposition}[section]
\newtheorem{remark}{Remark}

\newtheorem{ex}[theorem]{Example}

\newcommand{\MM}{\ensuremath{\mathcal M}}

\def\argmin{\mathop{\rm arg \; min}\limits}%
\newcommand{\ensnombre}[1]{\mathbb{#1}}
\newcommand{\N}{\ensnombre{N}}

\newcommand{\R}{\ensnombre{R}}

\renewcommand{\P}{\mathbb{P}}

\begin{document}
\title{New Fréchet features for random distributions and  associated sensitivity indices}
\author{Jean-Claude Fort$^{\rm a}$ and Thierry Klein $^{\rm b}$$^{\ast}$\thanks{ $^{\rm a}${MAP5 Universit\'e Paris Descartes, SPC, 45 rue des Saints P\`eres, 75006 Paris, France, $^{\rm b}${Institut de Math\'ematiques, University of Toulouse} $^\ast$Corresponding author. thierry.klein@math.univ-toulouse.fr}}}


\maketitle

\begin{abstract}
In this article we define new Fréchet features for random cumulative distribution functions using contrast. These contrasts allow to  construct Wasserstein costs and our new features minimize the average costs as the Fréchet mean minimizes the mean square Wasserstein$_2$ distance. 
An example of new features is the median, and more generally the quantiles. From these definitions, we are able to define sensitivity indices when the random distribution is the output of a stochastic code. Associated to the Fréchet mean we extend the Sobol indices, and in general  the indices associated to a contrast that we previously proposed.
\end{abstract}

\begin{bf}keywords\end{bf}
contrats, Wassertein costs, sensitivity index


\section*{Introduction}

Nowadays the output of many computer codes is not only a real multidimensional variable but frequently a function computed on so many points that it can be considered as a functional output. In particular this function may be the density or the cumulative distribution function ($c.d.f$) of a real random variable (phenomenon). In this article we focused on the case of a $c.d.f$ output. To analyze such outputs one  needs to choose a distance to compare various $c.d.f.$. Among the large possibilities offered by the literature we have chosen the Wasserstein distances (for more details on wasserstein distances we refer to \cite{villani2003tot}).  Actually for one dimensional probability distributions the Wasserstein$_p$ distance  simply is the $L^p$ distance of simulated random variables from a universal (uniform on $[0,1]$) simulator  $U$: $W_p^p(F,G)=\int_0^1|F^-(u)-G^-(u)|^p du=\mathbb E |F^-(U)-G^-(U)|^p$, where $F^-$ is the generalized inverse of $F$.  This means that using Wasserstein distances is to compare various $c.d.f$ from various codes on a same simulation space, which seems very natural in many situations. The most relevant cases seem to be $p=2$ and $p=1$, and in this paper we will work with.

In this article, we consider the problem of defining  a generalized notion of  barycenter of random probability measures on $\R$. It is a well known fact that the  set of  Radon probability measures endowed with the 2-Wasserstein distance is not an Euclidean space. Consequently, to define a notion of barycenter for  random probability measures, it is natural to use the notion of Fr\'echet mean \cite{fre} that is an extension of the usual Euclidean barycenter to non-linear spaces endowed with non-Euclidean metrics. If $\mathbb{Y}$ denotes a random variable with distribution $\P$  taking its value in a metric space $(\MM,d_{\MM})$, then a Fréchet mean (not necessarily unique) of the distribution $\P$  is a point $m^{\ast} \in \MM$ that is a global minimum (if any) of the functional 
$$
J(m) = \frac{1}{2} \int_{\MM} d^{2}_{\MM}(m,y) d\P(y) \quad \mbox{ i.e. }  \quad  m^{\ast} \in \argmin_{m \in \MM} J(m).
$$
In this paper, a Fréchet mean of a random variable $\mathbb{Y}$ with distribution $\P$ will be also called a barycenter.  For random variables belonging to nonlinear metric spaces, a well-known example is the computation of  the mean of a set of planar shapes in the Kendall's shape space \cite{kendall} that leads to the  Procrustean means studied in  \cite{MR1108330}. Many properties of the Fr\'echet mean in finite dimensional Riemannian manifolds (such as consistency and uniqueness) have been investigated in \cite{Afsari, MR2914758,batach1,batach2,BK2015,MR2816349}.   $ $\\

This article is an attempt  to use these tools and some extensions for analyzing computer codes outputs in a random environment, what is the subject of computer code experiments. In the first section we define new contrasts for  random $c.d.f.$ by considering  generalized " Wasserstein" costs. From this,  in the second section we define new features in the way of the Fr\'echet mean that we call Fr\'echet features. Then we propose some examples. The next two sections are devoted to a sensitivity analysis of random $c.d.f.$, first from a Sobol point of view that we generalized to a contrast point of view as in \cite{FKR2013}.

\section{Wasserstein distances and Wasserstein costs for unidimensional distributions}
For any $p\ge 1$ we may define a Wasserstein distance between two distribution of probability, denoted $F$ and $G$ (their cumulative distribution functions, $c.d.f.$) on $\mathbb R^d$ by: 

$$ W_p^p(F,G)=\min_{(X,Y)}\mathbb E\|X-Y\|^p,$$
where the random variables ($r.v.$'s) have $c.d.f.$ $F$ and $G$ ($X\sim F, Y\sim G$), assuming that $X$ and $Y$ have finite moments of order $p$. We call Wassertein$_p$ space the space of all $c.d.f.$ of $r.v.$'s with finite moments of order $p$.

As previously mentioned, in the unidimensional  case where $d=1$, it is well known that $W_p(F,G)$ is explicitly computed by:
$$ W_p^p(F,G)=\int_0^1|F^-(u)-G^-(u)|^p du=\mathbb E|F^-(U)-G^-(U)|^p.$$
Here $F^-$ and $G^-$ are the generalized inverses of $F$ and $G$ that are increasing  with limits $0$ and $1$, and $U$ is a $r.v.$ uniform on $[0,1]$. Of course $F^-(U)$ and $G^-(U)$ have $c.d.f.$ $F$ and $G$.

This result  extends 
to more general contrast functions. 
\begin{definition}
We call contrast functions any application $c$  from $\R^{2}$ to $\R$ satisfying the "measure  property" $\cal P$ defined by 
\begin{equation*}
\mathcal{ P}: \forall x\leq x'\ \mathrm{and\ } \forall y\leq y', c(x',y')-c(x',y)-c(x,y')+c(x,y)\leq 0,
\end{equation*}
 meaning that $c$ defines  a negative measure on $\mathbb R^2$.
\end{definition}
\begin{ex}
 $c(x,y)=-xy$ satisfies the $\cal P$ property.
\end{ex}
\begin{remark} If $c$ satisfies $\cal P$ then any function of the form $a(x)+b(y)+c(x,y)$ satisfies $\cal P$. For instance $(x-y)^2=x^2+y^2-2xy$ satisfies $\cal P$.

\end{remark}

\begin{remark} More generally if $C$ is a convex real function then c(x,y)=C(x-y) satisfies $\cal P$. This is the case of  $|x-y|^p$, $p\ge 1$.
\end{remark}

\begin{definition}
We define de Skorohod space $\mathcal{D}:=\mathcal{D}\left([0,1]\right)$ of all distribution functions that is the space of all non decreasing function from $\R$ to $[0,1]$ that are c\`ad-l\`ag with limit $0$ (resp. $1$) in $-\infty$ (resp. $+\infty$) equiped with the supremum norm.
\end{definition}

\begin{definition}[The $c-$Wasserstein cost]
For any $F\in\mathcal{D}$, any $G\in\mathcal{D}$  and any positive contrast  function $c$, we define the $c-$Wasserstein cost by 
$$\displaystyle W_c(F,G)=\min_{(X\sim F,Y\sim G)}\mathbb E \left(c(X,Y)\right)<+\infty$$

\end{definition}
The following theorem can be found in (\cite{Cambabis76}).

\begin{theorem}[Cambanis, Simon, Stout \cite{Cambabis76}] Let $c$ a function from $\mathbb R^2$ taking values in $\mathbb R$. Assume that it satisfies the "measure  property"  $\cal P $. Then

$$ W_c(F,G)=\int_0^1 c(F^-(u),G^-(u))du=\mathbb E\ c(F^-(U),G^-(U)),$$
where $U$ is a random variable uniformly distributed on $[0,1]$.
\end{theorem}

At this point we may notice that in a statistical framework one encounter many contrasts that are defined via a convex function. Actually many features of probability distribution can be characterized via such a contrast function. For instance an interesting case is the quantiles. Applying the previous remark we get:

\begin{proposition} For any $\alpha\in (0,1)$ the contrast function (pinball function) associated to the $\alpha$-quantile $c_{\alpha}(x,y)=(1-\alpha)(y-x) {\bold 1}_{x-y<0}+\alpha (x-y) {\bold 1}_{x-y\ge 0}$ satisfies $\cal P$.

\end{proposition}

This result is the starting point of the definition of some new features of random $c.d.f.$.

\section{Extension of the Fréchet mean to other features}

A Fréchet mean ${\cal E} X$ of a $r.v.$ $X$ taking values in a metric space $({\cal M}, d)$ is define as (whenever it exists):

$${\cal E} X\in \mbox{argmin}_{\theta \in {\cal M} } \mathbb E\ d(X,\theta)^2.$$

That means that it minimizes the contrast $\mathbb E\ d(X,\theta)^2$ which is an extension of the classical contrast $\mathbb E\|X-\theta\|^2$ in $\mathbb R^d$.

Adopting this point of view we can define a "Fréchet feature" associated to a convenient  contrast function. \\

Now we consider  a probability space $\left(\Omega,\mathcal{A},\P\right)$ and a measurable application  $\mathbb F$ from $\Omega$ to $\mathcal{D}$. Take $c$ a positive contrast (satisfying property $\mathcal{P}$) and define the analogously to the Fr\'echet mean, the  Fr\'echet feature associated to $c$ or contrasted by $c$ as it follows:

\
\begin{definition} Assume that $\mathbb F$  is a random  variable taking values in  $\mathcal{D}$. 
 Let $c$ be a non negative contrast function satisfying the property  $\cal P$. We define a $c$-contrasted feature ${\cal E}_c  \mathbb F$ of $\mathbb F$ by:

$${\cal E}_c  \mathbb F \in \mbox{argmin}_{G\in \mathcal{D}}\mathbb E\left(W_c(\mathbb F, G)\right).$$

\end{definition}

Of course this definition coincides with the Fréchet mean in the Wasserstein$_2$ space when using  the "contrast function" $c(F,G)=W_2^2(F,G)$.

\begin{theorem} \label {frechetfeature} If $c$ is a positive cost function satisfying  the property $\cal P$, if the application defined  on $(\omega,u)\in\Omega\times(0,1)$ by  $\mathbb F^{-}(\omega,u)$ is measurable  and if ${\cal E}_c  \mathbb F $ exists and is unique  we have:

$$ ({\cal E}_c  \mathbb F)^{-}(u) = \mbox{argmin}_{s\in \mathbb R} \mathbb E c(\mathbb F^{-}(u),s).$$

That is ${\cal E}_c \mathbb F$ is the inverse of the function taking value at $u$ the  $c$-contrasted feature of the real $r.v.$ $\mathbb F^{-}(u)$. For instance the Fréchet mean in the Wasserstein$_2$ space is the inverse of the function $u\longrightarrow \mathbb E\ \left( \mathbb F^{-}(u)\right)$.

\end{theorem}
\begin{remark}
Here, we proposed a general framework  on $\mathbb{F}$ and made some strong assumptions on existence uniqueness and measurability. But one can construct explicit parametric models for $\mathbb{F}$. We refer to \cite{BK2015} for such example. In particular in \cite{BK2015}, the authors used some results of \cite{MR2643592} that ensures measurability for some parametric models on $\mathbb{F}$.
\end{remark}
Another example is the Fréchet median. A contrast function defining the median in $\mathbb R$ is $ |x-y|$. An immediate extension to the Wassertein$_1$ space is to consider  the "contrast function" $c(F,G)=W_1(F,G)$. Thus we obtain the Fréchet median of a random $c.d.f.$ as :

$$(\mbox{Med} (\mathbb F))^{-}(u)\in \mbox{Med}(  \mathbb F^{-}(u)).$$

More generally we can define an $\alpha$-quantile of a random $c.d.f.$,  $q_\alpha(\mathbb F)$, as:

$$(q_\alpha (\mathbb F))^{-}(u)\in q_\alpha(  \mathbb F^{-}(u)),$$

where $q_\alpha(X)$ is the set of the $\alpha$-quantiles of $X$ taking its values in $\mathbb R$.\\

{\bf Proof} of Theorem \ref{frechetfeature}.\\

Since $c$ satisfies $\cal P$ we have: 

$$\mathbb E\  W_c(\mathbb F, G)=\mathbb E\int_0^1 c(\mathbb F^-(u),G^-(u))du=\int_0^1 \mathbb E\  c(\mathbb F^-(u),G^-(u))du,$$

by Fubini's theorem. 

Now for all $u\in (0,1)$ the quantity  $\mathbb E\ c(\mathbb F^-(u),G^-(u))$ is minimum for $G^-(u)$ a feature contrasted by $c$. Noticing that this results in an increasing and càd-làg function the theorem follows. \hfill $\square$

\section{Example}

In this section we illustrate our definitions through an example.\\ 

Let $F_0$ an increasing absolutely continuous  $c.d.f$ (hence $F_0^{-1}$ exists), $X$ a $r.v.$ with distribution $F_0$, $M$ and  $\Sigma$ two real  $r.v.$'s, $\Sigma$>0. We consider the random $c.d.f.$ $\mathbb F$ of $\Sigma X+M$. We have:  
$$\mathbb F(x)= F_0(\frac{x- M}{\Sigma}) \mbox{ and } \mathbb F^{-1}(u)=\Sigma F_0^{-1}(u)+M.$$

As well known the Fréchet mean of $\mathbb F$ is given by: $({\cal E} (\mathbb F))^{-1}(u)= \Sigma F_0^{-1}(u)+ M$, thus $\displaystyle{\cal E} (\mathbb F)(x)=F_0(\frac{x-\mathbb E M}{\mathbb E \Sigma})$.

Now using the $\alpha$-quantile contrast $c_{\alpha}(x,y)=(1-\alpha)(y-x) {\bold 1}_{x-y<0}+\alpha (x-y) {\bold 1}_{x-y\ge 0}$ and following our definition, we define the "$\alpha$-quantile" of $\mathbb F$: 

 $$(q_\alpha (\mathbb F))^{-1}(u)=q_\alpha(  \Sigma F_0^{-1}(u)+M).$$
 
Assuming that $\Sigma=1$ it simplifies in  $q_\alpha (\mathbb F)(x)=F_0(x-q_\alpha(M))$.  When $M=0$  we have $\displaystyle q_\alpha (\mathbb F)(x)=F_0(\frac{x}{q_\alpha(\Sigma)})$ (see figure(\ref{linear})).

Once these features defined, referring to computer experiment framework, in the next section we propose a sensitivity analysis of these Fr\'echet features of a random $c.d.f.$ as stochastic output of a computer code. 
 
 \section{Sensitivity indices for a random $c.d.f.$}
 
 \subsection{Sobol index}
 
 A very classical problem in the study of computer code experiments (see  \cite{sant:will:notz:2003}) is the evaluation of the relative influence of the input variables on some numerical result obtained by a computer code. 
This study is usually called sensitivity analysis in this paradigm and has been widely assessed (see for example \cite{sobol1993}, \cite{saltelli-sensitivity}, \cite{rocquigny2008uncertainty} and references therein).
More precisely,
the numerical result of interest $Y$ is seen as a function of the vector of the distributed input $(X_i)_{i=1,\cdots,d}$ ($d\in\N_*$). Statistically speaking, we are dealing here with the unnoisy non parametric model
\begin{equation}
Y=f(X_{1},\ldots, X_{d}),
\label{momodel}
\end{equation}
where $f$ is a regular unknown numerical function on the state space  $E_1\times E_2\times\ldots\times E_d$ on which the distributed variables $(X_{1},\ldots, X_{d})$ are living. Generally, the inputs are assumed to be stochastically independent and sensitivity analysis is performed by using the so-called Hoeffding decomposition (see \cite{van2000asymptotic} and \cite{anton84}). In this functional decomposition $f$ is expanded as a $L^2$ sum of uncorrelated functions involving only a part of the random inputs. For any subset $v$ of $I_d=\{1,\ldots,d\}$ this leads to an index called the Sobol index (\cite{sobol1993}) that measures the amount of {\it randomness} of $Y$ carried in the subset of input variables $(X_i)_{i\in v}$. Without loss of generality, let us consider the case where $v$ reduces to a singleton.
 Let us first recall some well known facts about Sobol index. The global Sobol index quantifies the influence of the $r.v.$ $X_i$ on 
the output $Y$. This index is based on the variance (see \cite{sobol1993},\cite{saltelli-sensitivity}): more precisely, it compares the total 
variance of $Y$ to the expected variance of the variable $Y$ conditioned by $X_i$,
\begin{equation}\label{sobol2}
S_i=\frac{\text{Var}(\mathbb E[Y|X_i])}{\text{Var}(Y)}.
\end{equation}
By the property of the conditional expectation it writes also 
\begin{equation}\label{sobol1}
S_i=\frac{\text{Var}(Y)-\mathbb E(\text{Var}[Y|X_i])}{\text{Var}(Y)}.
\end{equation}

In view of this formula we can define a Sobol index for the Fréchet mean of a random $c.d.f.$ $\mathbb F=h (X_1,\ldots,X_d)$. Actually we define $\text{Var}(\mathbb F)=\mathbb E W_2^2(\mathbb F,{\cal E}(\mathbb F))$, and 

$$S_i(F)=\frac{\text{Var}(\mathbb F)-\mathbb E(\text{Var}[\mathbb F|X_i])}{\text{Var} \mathbb F}.$$

From Theorem \ref{frechetfeature} we get:

 $$\text{Var}(\mathbb F)=\mathbb E\int_0^1 |\mathbb F^-(u)-{\cal E}(\mathbb F)^-(u)|^2du=\mathbb E\int_0^1 |\mathbb F^-(u)-{\mathbb E}\mathbb F^-(u)|^2du=\int_0^1\text{Var}(\mathbb F^-(u))du.$$
 
 And the Sobol index is now:
 
 $$S_i(\mathbb F)=\frac{\int_0^1\text{Var}(\mathbb F^-(u))du-\int_0^1\mathbb E \text{Var}[\mathbb F^-(u)|X_i]du}{\int_0^1\text{Var}(\mathbb F^-(u))du}=\frac{\int_0^1\text{Var}(\mathbb E[\mathbb F^-(u)|X_i])du}{{\int_0^1\text{Var}(\mathbb F^-(u))du}}.$$
 
 As a toy example, applying this to our previous example $\displaystyle \mathbb F(x)= F_0(\frac{x- M}{\Sigma})$, where $M$ and $\Sigma$ play the role of influent random variables, we find:
 
 $$S_{\Sigma}=\frac{{\text{Var }}\Sigma+2\text{cov}(\Sigma,M)\mathbb E \xi}{{\text{Var }}\Sigma+{\text{Var }}M+2\text{cov}(\Sigma,M) \mathbb E \xi },\\ S_{M}=\frac{{\text{Var }}M+2\text{cov}(\Sigma,M)\mathbb E \xi}{{\text{Var }}\Sigma+{\text{Var }}M+2\text{cov}(\Sigma,M) \mathbb E \xi }$$
 
 where $\xi$ has $c.d.f.$ $F_0$, since $\mathbb E \xi=\int_0^1 F_0^{-1}(u)du$.\\
 
 In practice $M$ and $\Sigma$ depends upon numerous random variables $(X_1,\ldots,X_d)$, then the Sobol index with respect to $X_i$ becomes:

 $$S_i=\frac{{\text{Var }}\mathbb E[\Sigma|X_i]+2\text{cov}(\mathbb E[\Sigma|X_i],\mathbb E[M|X_i])\mathbb E \xi+{\text{Var }}\mathbb E[M|X_i]}{{\text{Var }}\Sigma+{\text{Var }}M+2\text{cov}(\Sigma,M) \mathbb E \xi }$$

 \subsection{Sensitivity index associated to a contrast function}
 The formula (\ref{sobol1}) can be extended to more general contrast functions. The contrast function naturally associated to the mean of a real $r.v.$ is $c(x,y)=|x-y|^2$. We have $\mathbb E Y=\text{ argmin}_{\theta\in \mathbb R}\mathbb E c(Y,\theta)$ and $\text{Var}(Y)=\min_{\theta\in \mathbb R}\mathbb E c(Y,\theta)$. Thus the denominator of $S_i$ is the variation between the minimum value of the contrast and the expectation of the minimum of the same contrast when conditioning by the $r.v.$ $X_i$. Hence for a feature of  a real $r.v.$ associated to a contrast function $c$  we defined a sensitivity index (see (\cite{FKR2013})):
 
 $$S_{i,c}=\frac{\min_{\theta\in \mathbb R} \mathbb E c(Y,\theta)-\mathbb E \min_{\theta\in \mathbb R} \mathbb E [c(Y,\theta)|X_i]}{\min_{\theta\in\mathbb R} \mathbb E c(Y,\theta)}.$$

Along the same line, we now define a sensitivity index for a $c$-contrasted feature of a random $c.d.f.$ by:

$$S_{i,c}=\frac{\min_{G\in \mathbb W} \mathbb E W_c(\mathbb F,G)-\mathbb E \min_{G\in \mathbb W} \mathbb E [W_c(\mathbb F,G)|X_i]}{\min_{G\in\mathbb W} \mathbb E W_c(\mathbb F,G)}.$$

The computation of $S_{i,c}$ simplifies when $c$ satisfies the property $\cal P$ and assuming the uniqueness of ${\cal E}_c \mathbb F$:

$$S_{i,c}=\frac{\mathbb E\int_0^1c(\mathbb F^-(u),({\cal E}_c \mathbb F)^-(u))du  - \mathbb E[   \int_0^1c(\mathbb F^-(u),({\cal E}_c [\mathbb F|X_i])^-(u))du] }{\mathbb E\int_0^1c(\mathbb F^-(u),({\cal E}_c \mathbb F)^-(u)) }$$

where ${\cal E}_c [\mathbb F|X_i]$ is the c-contrasted feature conditional to $X_i$ ($i.e.$ with respect to the conditional distribution of $\mathbb F$), also assumed to be unique.

For instance if $c= |x-y|$, $({\cal E}_c \mathbb F)^-(u)$ is the "median" (assumed to be unique) of the random variable $\mathbb F^-(u)$ and:

$$S_{i,\text{Med}}=\frac{\mathbb E\int_0^1 |\mathbb F^-(u)-Med( \mathbb F^-(u))|du  - \mathbb E[   \int_0^1|\mathbb F^-(u)-Med [\mathbb F^-(u)|X_i]|du] }{\mathbb E\int_0^1 |\mathbb F^-(u)-Med( \mathbb F^-(u))|du }.$$

The same holds for any $\alpha$-quantile, using the corresponding contrast function $c_\alpha$ but whith less readable formula.

\section{Conclusion}
This article is an attempt to define interesting features for a functional output of a computer experiment, namely a random $c.d.f.$, together with its sensitivity analysis. This theory is based on contrast functions that allow to compute Wasserstein costs.  In the same way as the Fréchet mean for the Wassersstein$_2$ distance we have defined  features that minimize some contrasts made of these Wasserstein costs.
Straightforwardly  from the construction of that features we have developed a proposition of sensitivity analysis, first of Sobol type and then extended to sensitivity indices associated to our new  contrasts.
We intend to apply our methodology to an industrial problem: the PoD (Probability of Detection of a defect) in a random environment.  In particular we hope that our $\alpha$-quantiles will provide a relevantt tool to analyze that type of data.

\bibliographystyle{plain}
 \bibliography{Wcosts}
\end{document}